\documentclass{etds}
\newtheorem{theorem}{Theorem}

\newtheorem{cor}{Corollary}

\begin{document}
\ETDS{0}{0}{0}{0}

\runningheads{A.\ Kercheval}{Erratum for \itshape{Denjoy Minimal Sets
Are Far From Affine}}

\title{Erratum for \itshape{Denjoy Minimal Sets Are Far From Affine}}

\author{ALEC N.\ KERCHEVAL}

\address{Dept. of Mathematics, Florida State University,
Tallahassee, FL 32308-4510\\
\email{kercheva@math.fsu.edu}}

\recd{ April $2010$}

\begin{abstract} Theorem 2 of \cite{R1} is corrected by adding a $C^2$ bound to the hypotheses.
\end{abstract}
\medskip
Consider the circle {\bf R/Z} with coordinates $[0,1)$.
Let  $L \subset [0,1)$ be a compact interval and for $k \geq 2$ let ${\cal I} = \{I_1, \dots, I_k\}$ be a collection of pairwise
disjoint compact intervals with union $I \subset L$.
Define ${\cal S}^r({\cal I}, L)$ to be the set of $C^r$ functions $S: I \to L$ such that $|S'| >1$ on $I$, and
for each $j = 1, \dots, k$, $S[I_j] = L$. For such $S$, define its {\it nonlinearity} 
$
{\cal N}(S) \equiv \max_j \sup \{ \log (S'(x)/S'(y)): x,y \in I_j\}.
$

 Any $S \in {\cal S}^r({\cal I}, L)$ has a unique maximal invariant
(Cantor) set
$
C_S = \{x \in I : S^n(x) \in I \mbox{ for all } n > 0\}.
$
A Cantor set is {\it $C^1$-minimal}
if it is the minimal set of some $C^1$ diffeomorphism of the circle. The following
appears as Theorem 2 in \cite{R1}.
\begin{theorem}
Let $I_1, \dots, I_k, L$ be compact intervals as above.  Then there exists $\epsilon > 0$ (depending only
on $\{|I_j|/|L|: j = 1, \dots, k\}$) such that if $S \in {\cal S}^2({\cal I}, L)$ and ${\cal N}(S) < \epsilon$
then $C_S$ is not $C^1$-minimal.
\end{theorem}

A.\ Portela correctly points out that the proof in \cite{R1} only proves the following slightly weaker statement:
\begin{theorem} 
\label{new}
Let $I_1, \dots, I_k, L$ be compact intervals as above, and let $M>0$.  Then there exists $\epsilon > 0$ (depending only
on $\{|I_j|/|L|: j = 1, \dots, k\}$ and $M$) such that if $S \in {\cal S}^2({\cal I}, L)$, $|S''| < M$, and
$
{\cal N}(S)  < \epsilon
$
then $C_S$ is not $C^1$-minimal.
\end{theorem}

It still follows immediately that $C^1$-minimal Cantor sets are not $C^2$-nearly affine:

\begin{cor}
Let $I_1, \dots, I_k, L$ be given as above.  Then there exists $\epsilon >0$ such that
for all $A, S \in {\cal S}^r({\cal I}, L)$, if $A$ is locally affine and 
$||S - A||_{C^2} < \epsilon$, then $C_S$ is not $C^1$-minimal.
\end{cor}

\end{document}